\documentclass{amsart}
\usepackage{amsxtra,amssymb,multicol,graphicx}
\title{WHAT IS cost?}
\author{Damien Gaboriau}
\begin{document}

\rightline{\Huge \bf WHAT IS COST?}
\vskip1cm
\rightline{\it \huge Damien Gaboriau}
\vskip1in

\begin{multicols}{2}

  Orbit equivalence theory considers dynamical systems from the point of view of orbit equivalence relations.  The notion of {\em cost} is a useful invariant in this theory.

 When a countable group $\Gamma$ acts on a space it defines an equivalence relation: ``to be in the same orbit''. We consider measure preserving actions on a standard probability measure space. Cost was introduced by G. Levitt in order to quantify the amount of information needed to build this equivalence relation. 
 
% example of Z^2 action by rotations
Consider for instance the $\mathbb{Z}^2$-action on the circle 
$\mathbb{R}/2\pi \mathbb{Z}$ given by two rotations $a$ and $b$, whose angles together with $2\pi $ are rationally independent. Because of commutativity, there are many ways to check that two points are in the same orbit by using elementary {\em jumps} $x\sim a^{\pm 1} (x)$ and $x\sim b^{\pm 1}(x)$.
Indeed, the information encoded in the data $\{a,b\}$ is highly redundant.
Chooses instead some interval $I_\epsilon\subset \mathbb{R}/2\pi \mathbb{Z}$ of length $\epsilon>0$, and restricts the elementary jumps $x\sim b(x)$ to only those $x$'s in $I_\epsilon$ (and retain the $a$-jumps). Then $\Phi_{\epsilon}:=\{a, b{\vert I_\epsilon}\}$ still {\bf generates} the orbit equivalence relation $\mathcal{R}$ of the $\mathbb{Z}^2$-action: the smallest equivalence relation containing all the $(x,\varphi(x))$ for $\varphi\in \Phi_{\epsilon}$ and $x$ in the domain of $\varphi$, is $\mathcal R$ itself. 

In fact, because the $a$-orbits are dense, 
each point $z$ admits some $a$-iterate $a^n(z)$ in  $I_\epsilon$, so that the connection between $z$ and $b(z)$ may be recovered in $(2n+1)$-jumps, namely $n$ times $a$, followed by the restriction $b{\vert I_\epsilon}$ and then $n$ times $a^{-1}$.
The measures of the domains of $a$ and $b{\vert I_\epsilon}$ sum to $1+\epsilon$. This is by definition the {\bf cost} of $\Phi_\epsilon$.
And it is {\em cheaper} than the cost of $\{a,b\}$. Moreover, considering $\epsilon$ tending to $0$ leads one to declare $\mathcal{R}$ to have cost $=1$ (a priori $\leq 1$, but $\mathrm{Cost}(\mathcal{R})\geq 1$ when the classes are infinite).

% standard actions
More generally, consider an action $\alpha$ of a countable group $\Gamma$ on a standard Borel space $(X,\mu)$, preserving the probability measure $\mu$. Examples of such actions are plentiful, for instance Bernoulli shifts $(X,\mu)=(X_0,\mu_0)^{\Gamma}$ ($\Gamma$ acting by precomposition on functions $f:\Gamma\to X_0$ preserves the product measure $\mu_0^{\otimes\Gamma}$), or for instance the action by multiplication
of a countable subgroup of a compact group with its Haar measure.

In this measure theoretic context, all the constructions have to be measurable and sets of measure $0$ are neglected.
Assume that the  action is {\bf free}, i.e. the only element with a fixed point set of positive measure is the identity.

% def of cost
Consider a countable  family $\Phi=\{\varphi_j\}$ of isomorphisms $\varphi_j:A_j\to B_j$ between Borel subsets $A_j,B_j\subset X$ whose graphs are each contained in $\mathcal{R}_{\alpha}$, i.e. for each $j$, each $x\in A_j$ belongs to the $\alpha$-orbit of $\varphi_j(x)$.
The cost of $\Phi$ is the number of generators weighted by the measure of their domains, that is to say $\mathrm{Cost}(\Phi)=\sum_j \mu(A_j)$.

The {\bf cost} of the action $\alpha$, and equivalently of its orbit equivalence relation $\mathcal{R}_{\alpha}$, is defined as the infimum of the costs over all generating $\Phi$'s: 
\begin{center}$\mathrm{Cost}(\mathcal{R}_{\alpha}):=\inf\{\mathrm{Cost}(\Phi): \Phi \text{ generates } \mathcal{R}_{\alpha}\}.$\end{center}
It is clear, by taking the $\varphi_j=\alpha(\gamma_j)$ associated with a generating set $(\gamma_j)$ of $\Gamma$ that $\mathrm{Cost}(\mathcal{R}_{\alpha})$ is less than or equal to the {\bf rank} of $\Gamma$, i.e. its minimal number of generators.

% measure \nu
There is a pedantic way of defining the rank of a countable group: as the infimum of the measures of the generating subsets, for the natural measure on $\Gamma$, namely the counting-measure.
When applied to $\mathcal{R}_{\alpha}$, this gives an interesting interpretation of the cost. As a subset of $X\times X$, $\mathcal{R}_{\alpha}$ is simply the union of the graphs of the maps $\alpha(\gamma):X\to X$, for $\gamma\in \Gamma$.
There is a natural measure $\nu$ on $\mathcal{R}_{\alpha}$ defined by pushing forward $\mu$ by the maps $s^\gamma: x\mapsto (x, \gamma(x))$ and 
gathering together the various measures $s^\gamma_{*}\mu$.

The cost of $\mathcal{R}_{\alpha}$ is equivalently defined as the infimum of the $\nu$-measures of the subsets of $\mathcal{R}_{\alpha}$, that are not contained in any proper sub-equivalence relation.

% amenable
The $\mathbb{Z}^2$-action above, for which we computed the cost, belongs to a larger class of examples. Ornstein and Weiss proved the following remarkable result: for any  free action of an infinite {\em amenable} group (for instance commutative, or nilpotent, or solvable groups), the orbit equivalence relation may also be defined by {\em a  single} transformation  of the space $\psi:X\to X$. For every $\gamma\in \Gamma$ and almost every $x\in X$, there is a certain iterate $n(\gamma, x)$ such that $\psi^{n(\gamma, x)}(x)=\alpha(\gamma)(x)$. For all these actions, the infimum in the definition of the cost is in fact a minimum and equals $1$. This pair of properties in turn, when satisfied for some free action, implies amenability of $\Gamma$.
 
% cost 1 non treeable
The direct product of any two infinite groups $\Gamma_1\times \Gamma_2$, with at least one infinite-order element, is easily seen to produce only cost $1$ free actions 
(a straightforward elaboration of the above $\mathbb{Z}^2$-action example).
This allows one to produce plenty of actions whose cost is $1$ but for which the infimum is not attained for any generating $\Phi$: just consider  non-amenable direct products. Recall that containing the free group $\mathbf{F}_n$ on $n$ ($>1$) generators prevents a group from being amenable.

% treeings
When there exists a unique path of elementary jumps $x\sim \varphi^{\pm 1}(x)$ to connecting any two points in the same orbit,
$\Phi$ is called a {\bf treeing}. 
This notion is useful for computing the cost of some actions, and this is one of the main results in the theory: 
\textit{A treeing always realizes the cost of the equivalence relation it generates} [1].
If $\Phi $ is not a treeing, it is always possible to restrict some $\varphi\in \Phi$ to a Borel subset of its domain and nevertheless continue to generate the same equivalence relation. 
On the other hand, if $\Psi$ is a treeing, then restricting any $\psi\in \Psi$ to a subset of its domain breaks some connecting path and thus stops generating $\mathcal{R}_{\alpha}$.
A treeing is minimal in this sense. 
And the above statement claims that one cannot expect any 
other $\Phi$ to appear that will generate $\mathcal{R}_{\alpha}$ at a cheaper cost

% treed equiv. rel.
As a consequence, the cost of any free action of the free group $\mathbf{F}_n$ equals exactly $n$.
Indeed, the family of transformations associated with a free generating set of the group is a treeing. 
If two free actions of two free groups on $(X,\mu)$ are {\bf orbit equivalent} (i.e. define the same orbit equivalence relation), then the groups must have the same rank. The orbit equivalence relation {\em remembers} this rank. As another consequence, the cost of any free action of 
$\mathrm{PSL}(2,\mathbf{Z})$ is $1+\frac{1}{6}$. Recall that
$\mathrm{PSL}(2,\mathbf{Z})\simeq\mathbf{Z}/3\mathbf{Z}*\mathbf{Z}/2\mathbf{Z}$ and
$1+\frac{1}{6}=(1-\frac{1}{3})+(1-\frac{1}{2})$.

% no fundamental domain
Naive strategies to produce treeings usually collide with the following fundamental fact:
when an infinite countable group acts freely, there is no measurable way to pick one point in each orbit.
This is because such a set $D$ of selected points would have infinitely many pairwise disjoint translates $\gamma(D)$, all with the same measure, in a finite measure space: this is a impossible.

%Grushko
Grushko's theorem states that the rank of a free product equals the sum of the ranks of the factors.
Similarly, the cost of a free action of a free product $\Gamma_1*\Gamma_2$ equals the sum of the costs of the action restricted to the factors $\Gamma_i$. 

% Schreier
Recall Schreier's theorem: 
a finitely generated infinite normal subgroup of a free group $\mathbf{F}_n$ must have finite index. This theorem extends to those groups $\Gamma$ whose free actions have cost $>1$.
% compression formula
Moreover, Schreier's formula $(p-1)=i (n-1)$ relates the rank $n$ of the ambient free group $\mathbf{F}_n$ to the rank $p$ of a subgroup of finite index~$i$. This formula admits a counterpart in cost theory:
$[\mathrm{Cost}(\mathcal{R}_{\alpha}\vert A)-1] =  \mu(A)\ [\mathrm{Cost}(\mathcal{R}_{\alpha})-1]$ ({\bf compression formula}), where $\mathcal{R}_{\alpha}\vert A$ is the restriction  of the equivalence relation $\mathcal{R}_{\alpha}$ to some Borel subset $A\subset X$ that meets each orbit and that is equipped with the normalized restricted measure.

% open questions
Indeed, the computations of cost made thus far raise several open questions.
Cost seems to depend not on the particular free action but solely on the group.
Is this	 true in general ({\bf fixed price problem})?
There is moreover a strange coincidence with a numerical invariant, the first $\ell^2$ Betti number $\beta_1^{(2)}(\Gamma)$. Namely, it seems that $\beta_1^{(2)}(\Gamma)+1=\mathrm{Cost}(\mathcal{R}_{\alpha})$, when $\alpha$ is a free action of an infinite group $\Gamma$, although only the inequality $\leq$ has been proved.
In particular, is it true that actions of infinite Kazhdan property (T) groups have cost $=1$?

% conclusion
Orbit equivalence theory and cost are related to several other mathematical fields, like operator algebras, percolation on graphs, geometric group theory, descriptive set theory, etc. 
Much of the recent progress in von Neumann algebras and orbit equivalence was the result of a successful cross-pollination between these fields.

\def\refname{Further reading}

\end{multicols}

\end{document}